\documentclass[11pt]{amsart}

\usepackage[utf8]{inputenc}
\usepackage[T1]{fontenc}
\usepackage{lmodern}
\usepackage{amsmath,amssymb,amsthm,mathtools}
\usepackage{booktabs}
\usepackage[margin=1.15in]{geometry}
\usepackage{hyperref}

\newtheorem{theorem}{Theorem}[section]
\newtheorem{lemma}[theorem]{Lemma}
\newtheorem{proposition}[theorem]{Proposition}

\theoremstyle{definition}
\newtheorem{definition}[theorem]{Definition}
\newtheorem{question}[theorem]{Question}
\theoremstyle{remark}
\newtheorem{remark}[theorem]{Remark}

\theoremstyle{plain}
\newtheorem*{thmA}{Theorem A}
\newtheorem*{propB}{Proposition B}
\newtheorem*{propC}{Proposition C}
\newtheorem*{thmD}{Theorem D}
\newtheorem*{corE}{Corollary E}

\numberwithin{equation}{section}

\newcommand{\R}{\mathbb{R}}
\newcommand{\C}{\mathbb{C}}
\newcommand{\cp}{\operatorname{cap}}
\newcommand{\Per}{\operatorname{Per}}
\newcommand{\Res}{\operatorname*{Res}}

\begin{document}

\title[The Flucher--Rumpf problem for the interior Bernoulli problem]{On the Flucher--Rumpf isoperimetric question for the interior Bernoulli free boundary problem: a conformal reformulation, quantitative bounds and the second variation at the disc}

\author{Ababacar Sadikhe Djit\'e}
\address{Universit\'e Cheikh Anta Diop de Dakar, S\'en\'egal}
\email{ababacarsadikhe.djite@ucad.edu.sn}

\author{Diaraf Seck}
\address{Universit\'e Cheikh Anta Diop de Dakar, S\'en\'egal}
\email{diaraf.seck@ucad.edu.sn}

\date{November 2026}



\begin{abstract}
Let $c(\Omega)$ be the Bernoulli constant of a bounded domain $\Omega\subset\R^N$, that is, the infimum of the parameters $\lambda$ for which the interior Bernoulli free boundary problem $\Delta u=0$ in $\Omega\setminus K$, $u=0$ on $\partial\Omega$, $u=1$ and $|\nabla u|=\lambda$ on $\partial K$ admits a classical solution. Flucher and Rumpf asked whether $c(\Omega)\ge c(\Omega^*)$, where $\Omega^*$ is the ball of the same volume. The variational analogue was settled by Daners and Kawohl, but the variational constant differs from $c$ already for the ball, and the original question remains open. We prove the following.
(i) Through a conformal parametrisation of $\Omega\setminus K$ by an annulus, every classical solution in a simply connected planar domain with connected, filled free set satisfies $\lambda\ge\mu/R$, where $\pi R^2=|\Omega|$ and $\mu\approx 2.485268$; hence $c(\Omega)\ge 0.9142\,c(\Omega^*)$ for every bounded planar convex domain of class $C^1$.
(ii) In this class, the Flucher--Rumpf inequality is equivalent to an area inequality for analytic functions in an annulus, which also provides an explicit generator of exact solutions.
(iii) A Rellich--Pohozaev identity yields, in every dimension $N\ge2$ and for every strictly star-shaped domain, $c(\Omega)\ge\max\{\sqrt{\theta(\Omega)},\kappa_N\}\,c(\Omega^*)$, where $\theta(\Omega)\in(0,1]$ is an explicit geometric quantity equal to $1$ only for balls and $\kappa_N$ is an explicit constant with $\kappa_3\approx0.650$ and $\kappa_N\to1$ as $N\to\infty$.
(iv) We compute, at the level of a formal expansion, the second variation of $c$ at the disc along area-preserving deformations, $c(\Omega_\varepsilon)=e+A_k\varepsilon^2+O(\varepsilon^4)$, and show that $A_k>0$ for every $k\ge2$, with $A_1=0$.
(v) We show that two-point rearrangement (polarization) does not transport subsolutions of the problem, because it creates corners at which the gradient is unbounded.
All statements are supported by high-accuracy computations with the method of fundamental solutions.
\end{abstract}

\maketitle

%=====================================================================
\section{Introduction}
%=====================================================================

\subsection{The interior Bernoulli problem}
Let $\Omega\subset\R^N$, $N\ge2$, be a bounded domain and $\lambda>0$. The interior Bernoulli problem consists in finding a compact set $K\subset\Omega$ and a function $u$ such that
\begin{equation}\label{eq:bernoulli}
\Delta u=0 \ \text{in } \Omega\setminus K,\qquad u=0 \ \text{on } \partial\Omega,\qquad u=1 \ \text{and}\ |\nabla u|=\lambda \ \text{on } \partial K.
\end{equation}
A \emph{classical solution} is a pair $(K,u)$ such that $\partial K$ is of class $C^{1,\alpha}$, $u$ is continuous on $\overline{\Omega}\setminus K^\circ$, $u\in C^1(\mathcal N\setminus K^\circ)$ for a neighbourhood $\mathcal N$ of $\partial K$, and \eqref{eq:bernoulli} holds pointwise, the gradient on $\partial K$ being the limit from $\Omega\setminus K$. Problems of this type arise in ideal fluid flow, electrochemical machining and shape optimisation; see Flucher and Rumpf \cite{FR} and the survey \cite{Se16}.

For $K\subset\Omega$ compact we denote by
\[
\cp(K;\Omega)=\inf\Big\{\int_\Omega|\nabla v|^2\,dx:\ v\in W^{1,2}_0(\Omega),\ v\ge1 \text{ a.e. on } K\Big\}
\]
the condenser capacity of $K$ relative to $\Omega$, and by $u_K$ the capacitary potential, i.e.\ the minimiser. If $(K,u)$ is a classical solution, then $u=u_K$ on $\Omega\setminus K$ and, by Green's formula, $\cp(K;\Omega)=\int_{\partial K}|\nabla u|\,d\sigma=\lambda\Per(K)$.

For convex $\Omega$ there is a threshold: by the work of Beurling \cite{Be}, Acker \cite{Ac}, Henrot and Shahgholian \cite{HS1,HS2}, and Cardaliaguet and Tahraoui \cite{CT}, there exists $c(\Omega)>0$ such that \eqref{eq:bernoulli} has no solution for $\lambda<c(\Omega)$ and admits a classical solution with $K$ convex for every $\lambda\ge c(\Omega)$; moreover, the solution at $\lambda=c(\Omega)$ is unique \cite{CT}. In general we set
\begin{equation}\label{eq:defc}
c(\Omega)=\inf\{\lambda>0:\ \eqref{eq:bernoulli}\ \text{admits a classical solution}\},
\end{equation}
the \emph{Bernoulli constant} of $\Omega$ (with $c(\Omega)=+\infty$ if no classical solution exists). Since $u$ is invariant and $|\nabla u|$ scales like $t^{-1}$ under $x\mapsto tx$, one has $c(t\Omega)=c(\Omega)/t$, and $c$ is invariant under rigid motions.

\subsection{The question of Flucher and Rumpf}
In \cite[p.~202]{FR} Flucher and Rumpf asked whether
\begin{equation}\label{eq:FR}
c(\Omega)\ge c(\Omega^*),
\end{equation}
$\Omega^*$ being the ball with $|\Omega^*|=|\Omega|$. For the ball $B_R\subset\R^N$ the radial solutions are explicit. With $K=\overline{B_s}$ one has $u(r)=\log(R/r)/\log(R/s)$ if $N=2$ and $u(r)=(r^{2-N}-R^{2-N})/(s^{2-N}-R^{2-N})$ if $N\ge3$, so that on $\partial B_s$
\begin{equation}\label{eq:radial}
\lambda_{\rm rad}(s)=\frac{1}{s\log(R/s)}\ \ (N=2),\qquad
\lambda_{\rm rad}(s)=\frac{N-2}{s\,\big(1-(s/R)^{N-2}\big)}\ \ (N\ge3).
\end{equation}
Minimising in $s\in(0,R)$ gives
\begin{equation}\label{eq:cball}
c(B_R)=\frac{e}{R}\ \ (N=2,\ s=R/e),\qquad c(B_R)=\frac{(N-1)^{\frac{N-1}{N-2}}}{R}\ \ \big(N\ge3,\ s=R\,(N-1)^{-\frac1{N-2}}\big),
\end{equation}
in agreement with \cite{BS} (case $p=2$). For $\lambda>c(B_R)$ two radial solutions coexist (the elliptic and the hyperbolic one) and they merge at $\lambda=c(B_R)$: the Bernoulli constant is a turning point (fold) of the solution branch.

Several partial answers are known. Ly and Seck \cite{LS} proved, for the $p$-Laplacian and star-shaped domains, the existence of classical solutions for all parameters above an explicit threshold, which yields an explicit upper bound for the Bernoulli constant. Cardaliaguet and Tahraoui \cite{CT} gave an upper bound in terms of the harmonic radius and the lower bound $c(\Omega)\ge c(B_{\bar R})$, where $B_{\bar R}$ is the smallest ball containing $\Omega$; this lower bound degenerates for elongated domains. Bianchini and Salani \cite{BS} established a Brunn--Minkowski inequality for $c$ in the convex class, from which an Urysohn-type inequality follows: among convex sets of given mean width, the ball minimises $c$. As they observe, this does not imply \eqref{eq:FR}. Daners and Kawohl \cite{DK} considered the minimal parameter $\Lambda(\Omega)$ for which the Alt--Caffarelli functional
\[
J_\lambda(v)=\int_\Omega|\nabla v|^2\,dx+\lambda^2\,\big|\{v<1\}\big|,\qquad v\in W^{1,2}_0(\Omega),
\]
has a nontrivial minimiser, and proved by Schwarz symmetrisation that $\Lambda(\Omega)\ge\Lambda(\Omega^*)$, with equality essentially only for balls. However, $\Lambda$ and $c$ are different quantities: for the disc $\Lambda(B_R)=2\sqrt e/R\approx3.30/R$, whereas $c(B_R)=e/R\approx2.72/R$; for $c(B_R)\le\lambda<\Lambda(B_R)$ solutions exist, but none of them has lower energy than $v\equiv0$. Inequality \eqref{eq:FR} itself remains open, even for convex planar domains; see also \cite{BS,CT}.

\subsection{Main results}
In the plane we use conformal mapping. Throughout Sections~\ref{sec:conformal}--\ref{sec:second}, $\Omega\subset\C$ is a bounded simply connected domain, and we only consider solutions whose free set $K$ is connected and filled ($\C\setminus K$ connected); we call them \emph{admissible solutions}. For convex $\Omega$ of class $C^1$, the solution at $\lambda=c(\Omega)$ is admissible by \cite{HS2}.

\begin{thmA}[quantitative isoperimetric bound in the plane]
Let $\Omega\subset\C$ be bounded and simply connected with $|\Omega|=\pi R^2$, and let $(K,u)$ be an admissible classical solution of \eqref{eq:bernoulli}. Then
\[
\lambda\ \ge\ \frac{\mu}{R},\qquad \mu=\min_{0<\rho<1}\frac{\sqrt{1-\rho^2}}{\rho\log(1/\rho)}=\frac{1}{\rho_*\sqrt{1-\rho_*^2}}\approx 2.485268,
\]
where $\rho_*\approx0.450764$ is the unique root in $(0,1)$ of $\log(1/\rho)=1-\rho^2$. In particular, for every bounded convex domain $\Omega\subset\R^2$ of class $C^1$,
\[
c(\Omega)\ \ge\ \frac{\mu}{e}\,c(\Omega^*)\ \approx\ 0.914279\,c(\Omega^*).
\]
\end{thmA}

To state the converse we introduce, for $0<\rho<1$, the annulus $A_\rho=\{\rho<|w|<1\}$ and the following class.

\begin{definition}\label{def:class}
Let $\mathcal A$ be the set of triples $(\rho,m,\varphi)$ with $\rho\in(0,1)$, $m>0$ and $\varphi$ analytic in $A_\rho$, continuous on $A_\rho\cup\{|w|=\rho\}$, with $\operatorname{Re}\varphi=0$ on $\{|w|=\rho\}$, such that $m e^{\varphi}$ has a primitive $F$ on $A_\rho$ (equivalently, $\Res_{w=0}e^\varphi=0$ in the Laurent sense) which is bounded and extends to an injective continuous map on $\{\rho\le|w|<1\}$. We write $K_F$ for the compact set bounded by the Jordan curve $F(\{|w|=\rho\})$ and $\Omega_F=F(A_\rho)\cup K_F$.
\end{definition}

\begin{propB}[equivalent extremal problem]
The Flucher--Rumpf inequality $\lambda\ge e/R$ holds for all admissible solutions in all bounded simply connected planar domains if and only if, for every $(\rho,m,\varphi)\in\mathcal A$,
\begin{equation}\label{eq:extremal}
|\Omega_F|\ \ge\ \pi m^2\,\big(e\rho\log(1/\rho)\big)^2 .
\end{equation}
Moreover, every element of $\mathcal A$ produces an explicit admissible solution in $\Omega_F$.
\end{propB}

Proposition~B is a reformulation rather than a new inequality: all the difficulty of \eqref{eq:FR} is transferred to the univalence constraint in Definition~\ref{def:class}, which is global and non-convex. Its interest lies in the converse part, which generates exact solutions on which \eqref{eq:FR} can be tested (Section~\ref{sec:numerics}).

\begin{propC}[second variation at the disc; formal]
Let $k\ge1$ and $\Omega_\varepsilon=\{r<(1+\varepsilon\cos k\theta)/\sqrt{1+\varepsilon^2/2}\}$, so that $|\Omega_\varepsilon|=\pi$. At the level of a formal expansion (see Remark~\ref{rem:rigour}),
\[
c(\Omega_\varepsilon)=e+A_k\varepsilon^2+O(\varepsilon^4),\qquad
A_k=\frac{e\,(k^2-1)\big[(k-1)-2(k^2-1)x_k-(k+1)x_k^2\big]}{2\big[(k-1)+(k+1)x_k\big]^2},
\]
where $x_k=e^{-2k}$.
One has $A_1=0$ (translations), $A_k>0$ for all $k\ge2$, and $A_k\sim e(k+1)/2$ as $k\to\infty$.
\end{propC}

The expansion is even in $\varepsilon$, because the rotation $\theta\mapsto\theta+\pi/k$ maps $\Omega_\varepsilon$ onto $\Omega_{-\varepsilon}$. Proposition~C is derived by two independent computations (the conformal parametrisation and a domain expansion in the physical plane), which agree identically and match the numerical values to four digits.

Our last main result is not restricted to the plane. For a bounded domain $\Omega\subset\R^N$ of class $C^{1,\alpha}$ which is strictly star-shaped with respect to the origin, i.e.\ $x\cdot\nu(x)>0$ on $\partial\Omega$, set
\[
I(\Omega)=\int_{\partial\Omega}\frac{d\sigma}{x\cdot\nu},\qquad \theta(\Omega)=\frac{I(\Omega^*)}{I(\Omega)},
\]
where $\Omega^*$ is the ball centred at the origin with $|\Omega^*|=|\Omega|$. We show in Section~\ref{sec:rellich} that $\theta(\Omega)\in(0,1]$, with $\theta(\Omega)=1$ if and only if $\Omega$ is a ball centred at the origin. Define also $\kappa_2=0$ and, for $N\ge3$,
\begin{equation}\label{eq:kappa}
\kappa_N=\sqrt{\frac{N-2}{N}}\;\frac{N\,(N/2)^{\frac1{N-2}}}{(N-1)^{\frac{N-1}{N-2}}}.
\end{equation}

\begin{thmD}[Rellich--Pohozaev bounds in every dimension]
Let $N\ge2$, let $\Omega\subset\R^N$ be a bounded domain of class $C^{1,\alpha}$, strictly star-shaped with respect to the origin, and let $(K,u)$ be a classical solution of \eqref{eq:bernoulli} with parameter $\lambda$ (no assumption is made on the topology of $K$). Then
\begin{equation}\label{eq:rellich}
\lambda^2\,N|K|=\int_{\partial\Omega}|\nabla u|^2\,(x\cdot\nu)\,d\sigma+(N-2)\cp(K;\Omega),
\end{equation}
and
\[
\lambda\ \ge\ \max\big\{\sqrt{\theta(\Omega)},\ \kappa_N\big\}\;c(\Omega^*).
\]
Consequently $c(\Omega)\ge\max\{\sqrt{\theta(\Omega)},\kappa_N\}\,c(\Omega^*)$.
\end{thmD}

\begin{table}[ht]
\centering
\begin{tabular}{@{}lcccccc@{}}
\toprule
$N$ & 3 & 4 & 5 & 10 & 20 & 100\\
\midrule
$\kappa_N$ & 0.6495 & 0.7698 & 0.8278 & 0.9234 & 0.9636 & 0.9930\\
\bottomrule
\end{tabular}
\caption{The constant $\kappa_N$ of Theorem~D; $\kappa_N\to1$ as $N\to\infty$.}
\label{tab:kappa}
\end{table}

For $N\ge3$, Theorem~D gives a uniform bound of the type of Theorem~A for all strictly star-shaped domains, convex or not. For $N=2$ the coefficient $N-2$ in \eqref{eq:rellich} vanishes; this is why a uniform constant in the plane requires a different tool, here the conformal parametrisation. In the plane, the two results are complementary.

\begin{corE}
Let $\Omega\subset\R^2$ be a bounded convex domain of class $C^{1,\alpha}$ and $x_0\in\Omega$. Then
\[
c(\Omega)\ \ge\ \max\Big\{\frac{\mu}{e},\ \sqrt{\theta_{x_0}(\Omega)}\Big\}\,c(\Omega^*),
\]
where $\theta_{x_0}$ is computed with respect to $x_0$. For an ellipse of aspect ratio $q$, centred at $x_0$, $\theta=2/(q+q^{-1})$, and the second bound is the better one for $q<1.8529$.
\end{corE}

Finally, Section~\ref{sec:polarization} explains why the natural symmetrisation route through polarization fails at the level of subsolutions, while it does work at the level of energies.

\subsection*{Organisation}
Section~\ref{sec:conformal} contains the conformal dictionary and the proofs of Theorem~A and Proposition~B. Section~\ref{sec:second} is devoted to Proposition~C. Section~\ref{sec:rellich} proves Theorem~D and Corollary~E. Section~\ref{sec:polarization} deals with polarization. Section~\ref{sec:numerics} reports the numerical experiments, and Section~\ref{sec:open} lists open problems.

%=====================================================================
\section{The conformal dictionary and the proofs of Theorem A and Proposition B}\label{sec:conformal}
%=====================================================================

\subsection{From solutions to analytic functions}
Let $(K,u)$ be an admissible classical solution in a bounded simply connected domain $\Omega\subset\C$. Since $K$ is connected and filled and $\partial K$ is a $C^{1,\alpha}$ curve, $K$ has nonempty interior and $\partial K$ is a single $C^{1,\alpha}$ Jordan curve. As $\partial\Omega$ is a nondegenerate continuum, $\Omega\setminus K$ is a doubly connected domain with two nondegenerate boundary components. Hence there exist a unique $\rho\in(0,1)$ and a conformal map
\[
F:\ A_\rho=\{\rho<|w|<1\}\longrightarrow\Omega\setminus K
\]
such that $F(w)\to\partial\Omega$ as $|w|\to1$ and $F(w)\to\partial K$ as $|w|\to\rho$.

\begin{lemma}[boundary regularity at the inner circle]\label{lem:boundary}
$F$ extends to an injective map of class $C^{1,\alpha}$ on $\{\rho\le|w|<1\}$, and $F'$ does not vanish on $\{|w|=\rho\}$.
\end{lemma}

\begin{proof}
The Kellogg--Warschawski theorem \cite[Thm.~3.6]{Po} concerns simply connected domains, so we reduce to that case. Fix $z_0\in K^\circ$ and let $T(z)=1/(z-z_0)$. Then $W=T(\widehat\C\setminus K)$ is a bounded simply connected domain containing $0=T(\infty)$, whose boundary $T(\partial K)$ is a $C^{1,\alpha}$ Jordan curve. Let $g:\mathbb D\to W$ be a Riemann map. By the Kellogg--Warschawski theorem, $g\in C^{1,\alpha}(\overline{\mathbb D})$ and $g'\neq0$ on $\overline{\mathbb D}$. Moreover $T(\Omega\setminus K)=W\setminus E$, where $E=T(\widehat\C\setminus\Omega)$ is a compact connected subset of $W$.

The map $h=g^{-1}\circ T\circ F$ is a conformal map of $A_\rho$ onto the doubly connected domain $\mathbb D\setminus g^{-1}(E)$, and $|h(w)|\to1$ as $|w|\to\rho$. By the Schwarz reflection principle, $h$ extends analytically to a neighbourhood of $\{|w|=\rho\}$ by $h(w)=1/\overline{h(\rho^2/\bar w)}$. This extension is injective near the circle: the two half-collars are mapped into $\mathbb D$ and into $\widehat\C\setminus\overline{\mathbb D}$ respectively, and the boundary map is a homeomorphism of circles. Hence $h'\neq0$ on $\{|w|=\rho\}$.

Near the circle, $T\circ F=g\circ h$ stays away from $0$, so $F=z_0+1/(g\circ h)$ is of class $C^{1,\alpha}$ up to $\{|w|=\rho\}$. Finally, $F'=-(g'\circ h)\,h'/(g\circ h)^2\neq0$ there.
\end{proof}

The function $w\mapsto\log|w|/\log\rho$ is harmonic in $A_\rho$, equal to $0$ on $|w|=1$ and to $1$ on $|w|=\rho$. The function $u\circ F$ is bounded and harmonic in $A_\rho$, with boundary limits $0$ as $|w|\to1$ (since $u$ is continuous up to $\partial\Omega$) and $1$ as $|w|\to\rho$. By uniqueness for the Dirichlet problem,
\begin{equation}\label{eq:uF}
u\circ F(w)=\frac{\log|w|}{\log\rho},\qquad w\in A_\rho.
\end{equation}
Set $L=\log(1/\rho)$. By conformal invariance, $|\nabla u(F(w))|\,|F'(w)|=|\nabla(u\circ F)(w)|=1/(|w|L)$. Letting $|w|\to\rho$ and using Lemma~\ref{lem:boundary}, we obtain on $\partial K$
\begin{equation}\label{eq:gradF}
|\nabla u(F(w))|=\frac{1}{\rho L\,|F'(w)|},\qquad |w|=\rho.
\end{equation}
The Bernoulli condition $|\nabla u|=\lambda$ is therefore equivalent to
\begin{equation}\label{eq:modulus}
|F'(w)|\equiv m\ \ \text{on } |w|=\rho,\qquad \lambda=\frac{1}{m\rho L}.
\end{equation}

\begin{lemma}\label{lem:mean}
Let $h=\log|F'|$, which is harmonic in $A_\rho$. Then for every $r\in(\rho,1)$
\[
\frac1{2\pi}\int_0^{2\pi}h(re^{i\theta})\,d\theta=\log m .
\]
\end{lemma}

\begin{proof}
For a harmonic function in an annulus, the circular mean $M(r)$ satisfies $(rM'(r))'=0$, so $M(r)=a+\beta\log r$ with
\[
\beta=rM'(r)=\frac1{2\pi}\int_0^{2\pi}r\,\partial_rh(re^{i\theta})\,d\theta=\frac1{2\pi}\oint_{|w|=r}d\arg F'(w),
\]
by the Cauchy--Riemann equations ($r\partial_r\log|f|=\partial_\theta\arg f$). Thus $\beta$ is the winding number of $F'$ along $|w|=r$.

The curve $\gamma_r=F(\{|w|=r\})$ is an analytic Jordan curve which separates $K$ from $\partial\Omega$. Since $F$ preserves orientation and maps the bounded region $\{\rho<|w|<r\}$ onto the region between $\partial K$ and $\gamma_r$, the curve $\gamma_r$ is positively oriented. Its tangent vector is $\partial_\theta F(re^{i\theta})=ire^{i\theta}F'(re^{i\theta})$. By the Umlaufsatz, the argument of the tangent increases by $2\pi$ along $\gamma_r$. The argument of $ire^{i\theta}$ also increases by $2\pi$, hence the winding number of $F'$ is $0$.

Thus $M$ is constant on $(\rho,1)$. By Lemma~\ref{lem:boundary}, $h$ is uniformly continuous on $\{\rho\le|w|\le r_1\}$ for every $r_1<1$, so $M(r)\to\frac1{2\pi}\int_0^{2\pi}h(\rho e^{i\theta})\,d\theta=\log m$ as $r\to\rho$.
\end{proof}

Lemma~\ref{lem:mean} shows that $\log F'$ has a single-valued branch in $A_\rho$. Writing $F'=me^{\varphi}$, the function $\varphi$ is analytic in $A_\rho$, continuous up to $|w|=\rho$, and $\operatorname{Re}\varphi=0$ there. In Laurent series this means
\begin{equation}\label{eq:laurent}
\varphi(w)=i\gamma+\sum_{n\ge1}\big(c_nw^n-\overline{c_n}\,\rho^{2n}w^{-n}\big),\qquad\gamma\in\R,
\end{equation}
since on $|w|=\rho$ one has $w^{-n}=\bar w^{\,n}\rho^{-2n}$ and each bracket equals $2i\operatorname{Im}(c_nw^n)$. Single-valuedness of the primitive $F$ is the condition $\Res_{w=0}e^\varphi=0$. Since $F$ is bounded and, by Lemma~\ref{lem:boundary}, injective on $\{\rho\le|w|<1\}$, we have $(\rho,m,\varphi)\in\mathcal A$, $K=K_F$ and $\Omega=\Omega_F$.

\subsection{From analytic functions to solutions}\label{sec:converse}
Conversely, let $(\rho,m,\varphi)\in\mathcal A$ with primitive $F$. Since $\operatorname{Re}\varphi=0$ on $|w|=\rho$, the Schwarz reflection $\varphi(w)=-\overline{\varphi(\rho^2/\bar w)}$ extends $\varphi$ analytically across $\{|w|=\rho\}$. Hence $F$ is analytic across that circle with $|F'|=m\neq0$ there, and $\Gamma_\rho=F(\{|w|=\rho\})$ is an analytic Jordan curve.

Let $K_F$ be the closed region it bounds and $\Omega_F=F(A_\rho)\cup K_F$. Then $\Omega_F$ is open, bounded (as $F$ is bounded) and simply connected: its complement is the unbounded component of $\C\setminus F(A_\rho)$. Define $u=\log|F^{-1}|/\log\rho$ in $F(A_\rho)$ and $u=1$ on $K_F$. Then $u$ is harmonic in $\Omega_F\setminus K_F$ and of class $C^1$ up to $\Gamma_\rho$, since $F^{-1}$ is analytic up to $\Gamma_\rho$. As $z\to\partial\Omega_F$, $F^{-1}(z)$ leaves every compact subset of $\{\rho\le|w|<1\}$, so $|F^{-1}(z)|\to1$ and $u(z)\to0$. By \eqref{eq:gradF}--\eqref{eq:modulus}, $(K_F,u)$ is an admissible classical solution of \eqref{eq:bernoulli} in $\Omega_F$ with $\lambda=1/(m\rho L)$. Moreover
\begin{equation}\label{eq:area}
|\Omega_F|=|K_F|+\iint_{A_\rho}|F'(w)|^2\,dA(w).
\end{equation}
This gives an explicit way to generate exact solutions, which we use in Section~\ref{sec:numerics}.

\subsection{Proof of Theorem A}
By Jensen's inequality and Lemma~\ref{lem:mean}, for every $r\in(\rho,1)$,
\[
\frac1{2\pi}\int_0^{2\pi}|F'(re^{i\theta})|^2\,d\theta=\frac1{2\pi}\int_0^{2\pi}e^{2h(re^{i\theta})}\,d\theta\ \ge\ \exp\Big(\frac1{\pi}\int_0^{2\pi}h(re^{i\theta})\,d\theta\Big)=m^2 .
\]
Integrating in $r$ and using \eqref{eq:area} with $|K|\ge0$,
\[
\pi R^2=|\Omega|\ \ge\ |\Omega\setminus K|=\int_\rho^1\!\!\int_0^{2\pi}|F'(re^{i\theta})|^2\,r\,d\theta\,dr\ \ge\ \pi m^2(1-\rho^2).
\]
Hence $m\le R/\sqrt{1-\rho^2}$ and, by \eqref{eq:modulus},
\[
\lambda=\frac1{m\rho L}\ \ge\ \frac{\sqrt{1-\rho^2}}{R\,\rho\log(1/\rho)}\ \ge\ \frac{\mu}{R}.
\]
It remains to compute $\mu$. Let $f(\rho)=\sqrt{1-\rho^2}/(\rho\log(1/\rho))$, which tends to $+\infty$ as $\rho\to0$ and as $\rho\to1$. A direct computation gives
\[
\frac{f'}{f}(\rho)=-\frac{\rho}{1-\rho^2}-\frac1\rho+\frac{1}{\rho\log(1/\rho)},
\]
which vanishes if and only if $\log(1/\rho)=1-\rho^2$. The function $\psi(\rho)=\log(1/\rho)-1+\rho^2$ satisfies $\psi'(\rho)=2\rho-1/\rho$. It is therefore decreasing on $(0,2^{-1/2})$ and increasing on $(2^{-1/2},1)$; it tends to $+\infty$ at $0$, is negative at $2^{-1/2}$ and vanishes at $\rho=1$. Hence it has a unique root $\rho_*\in(0,2^{-1/2})$, and at this root $f(\rho_*)=1/(\rho_*\sqrt{1-\rho_*^2})$. Numerically, $\rho_*=0.4507636\ldots$ and $\mu=2.4852675\ldots$.

If $\Omega$ is convex of class $C^1$, then by \cite{HS2} problem \eqref{eq:bernoulli} has a classical solution with $K$ convex at $\lambda=c(\Omega)$. This solution is admissible, so $c(\Omega)\ge\mu/R=(\mu/e)\,c(\Omega^*)$ by \eqref{eq:cball}. \qed

\begin{remark}\label{rem:lossy}
Two inequalities are lossy in this argument: Jensen's inequality, which is an equality only if $|F'|$ is constant on each circle, and the bound $|K|\ge0$. For the disc ($F(w)=mw$, $\rho=1/e$) one has $|K|=\pi m^2\rho^2$. Any lower bound of the form $|K|\ge\pi m^2\rho^2\,\Theta$, with $\Theta>0$ depending on how far $\frac1{2\pi}\int|F'|^2$ is from $m^2$, would improve the constant $\mu/e$.
\end{remark}

\subsection{Proof of Proposition B}
Assume that \eqref{eq:extremal} holds on $\mathcal A$, and let $(K,u)$ be an admissible solution in $\Omega$ with $|\Omega|=\pi R^2$. The dictionary above gives $(\rho,m,\varphi)\in\mathcal A$ with $\Omega=\Omega_F$ and $\lambda=1/(m\rho L)$. Then \eqref{eq:extremal} reads $\pi R^2\ge\pi(e/\lambda)^2$, i.e.\ $\lambda\ge e/R$.

Conversely, assume the Flucher--Rumpf inequality for admissible solutions. Each $(\rho,m,\varphi)\in\mathcal A$ produces, by Section~\ref{sec:converse}, the admissible solution $(K_F,u)$ in $\Omega_F$ with $\lambda=1/(m\rho L)$. Hence $\lambda\ge e/R$ with $\pi R^2=|\Omega_F|$, that is, $|\Omega_F|\ge\pi e^2/\lambda^2=\pi m^2(e\rho L)^2$. \qed

\begin{remark}\label{rem:ezero}
The disc corresponds to $\varphi\equiv i\gamma$ and $\rho=1/e$, where \eqref{eq:extremal} is an equality. Inequality \eqref{eq:extremal} with the factor $(e\rho\log(1/\rho))^2$ replaced by $1$ is false for $\rho$ close to $1$ (see \eqref{eq:Q} below, where $Q_k(\rho)<0$ for $\rho$ near $1$). Thus the factor $(e\rho\log(1/\rho))^2\le1$, which vanishes as $\rho\to1$, is essential.
\end{remark}

%=====================================================================
\section{Second variation at the disc}\label{sec:second}
%=====================================================================

\subsection{Reduction at the turning point}
Let $\Omega_\varepsilon$ be as in Proposition~C. Along the branch of solutions in $\Omega_\varepsilon$, parametrised by a smooth parameter $s$ with $s=s_0$ at the turning point of the disc, we write
\[
\lambda(s,\varepsilon)=\lambda_0(s)+\varepsilon\lambda_1(s)+\varepsilon^2\lambda_2(s)+O(\varepsilon^3).
\]
The rotation $\theta\mapsto\theta+\pi/k$ maps $\Omega_\varepsilon$ onto $\Omega_{-\varepsilon}$, so $\lambda_1\equiv0$ and the expansion is even in $\varepsilon$. Since $\lambda_0'(s_0)=0$ and $\lambda_0''(s_0)>0$, the minimum in $s$ moves by $O(\varepsilon^2)$, which changes the minimal value only by $O(\varepsilon^4)$. Hence
\begin{equation}\label{eq:fold}
c(\Omega_\varepsilon)=\min_s\lambda(s,\varepsilon)=e+\varepsilon^2\lambda_2(s_0)+O(\varepsilon^4).
\end{equation}
For the same reason $\lambda_2(s_0)$ does not depend on the choice of parametrisation of the branch. This is the simplest instance of a Lyapunov--Schmidt reduction at a fold.

\subsection{First derivation: conformal parametrisation}
We use the conformal modulus $\rho$ as the branch parameter ($\rho_0=1/e$). Take $k\ge2$, $a>0$ and
\[
\varphi_\varepsilon(w)=\varepsilon a\,(w^k-\rho^{2k}w^{-k}),
\]
which is of the form \eqref{eq:laurent}. All exponents in $e^{\varphi_\varepsilon}=\sum_jg_jw^j$ are multiples of $k$, so $g_{-1}=0$ and
\[
F(w)=m\sum_j\frac{g_j}{j+1}\,w^{j+1}
\]
is single-valued and analytic in $\C\setminus\{0\}$. For $\varepsilon$ small it is a $C^1$-small perturbation of $mw$ on a neighbourhood of $\overline{A_\rho}$, hence univalent there. Thus $(\rho,m,\varphi_\varepsilon)\in\mathcal A$ and, by Section~\ref{sec:converse}, we obtain an exact solution with $\lambda=1/(m\rho L)$.

\emph{Area.} If $F(e^{i\theta})=\sum_na_ne^{in\theta}$, the area enclosed by the Jordan curve $F(\{|w|=1\})$ is $|\Omega_F|=\pi\sum_nn|a_n|^2$. Here $a_{j+1}=mg_j/(j+1)$. From $e^{\varphi}=1+\varphi+\varphi^2/2+O(\varepsilon^3)$ we get
\[
g_0=1-\varepsilon^2a^2\rho^{2k}+O(\varepsilon^4),\quad g_k=\varepsilon a+O(\varepsilon^3),\quad g_{-k}=-\varepsilon a\rho^{2k}+O(\varepsilon^3),\quad g_{\pm2k}=O(\varepsilon^2),
\]
and therefore
\begin{equation}\label{eq:Q}
\frac{|\Omega_F|}{\pi m^2}=\sum_{j\ne-1}\frac{|g_j|^2}{j+1}=1+\varepsilon^2a^2Q_k(\rho)+O(\varepsilon^4),\qquad Q_k(\rho)=\frac1{k+1}-2\rho^{2k}-\frac{\rho^{4k}}{k-1}.
\end{equation}

\emph{Shape of the outer boundary.} On the outer circle,
\[
F(e^{i\theta})=m\,e^{i\theta}\Big(1+\varepsilon a\Big[\frac{e^{ik\theta}}{k+1}+\frac{\rho^{2k}e^{-ik\theta}}{k-1}\Big]+O(\varepsilon^2)\Big).
\]
Up to an $O(\varepsilon)$ tangential reparametrisation, $\partial\Omega_F$ is therefore the curve $r=m\,(1+\delta\cos k\psi+O(\varepsilon^2))$, with
\[
\delta=\varepsilon a\,S_k(\rho),\qquad S_k(\rho)=\frac1{k+1}+\frac{\rho^{2k}}{k-1}.
\]

\emph{Matching.} The $O(\varepsilon^2)$ corrections to $\varphi$ needed to match a prescribed outer boundary exactly involve only the modes $0$ (with imaginary coefficient, i.e.\ a rotation) and $\pm2k$. They do not affect \eqref{eq:Q} at order $\varepsilon^2$. A mode $2k$ of amplitude $O(\delta^2)$ in the boundary changes $c$ only at order $\delta^4$: by rotational invariance, the cubic terms of the expansion of $c$ couple it only to the square of the $k$-th mode. Normalising $|\Omega_F|=\pi$ and using \eqref{eq:modulus}, the solution of conformal modulus $\rho$ in the domain with boundary amplitude $\delta$ satisfies
\[
\lambda(\rho,\delta)^2=\frac{1+\delta^2Q_k(\rho)/S_k(\rho)^2}{(\rho\log(1/\rho))^2}+O(\delta^4).
\]
By \eqref{eq:fold}, $c=e\big(1+\delta^2Q_k(1/e)/(2S_k(1/e)^2)\big)+O(\delta^4)$, that is,
\begin{equation}\label{eq:Ak1}
A_k=\frac{e\,Q_k(1/e)}{2\,S_k(1/e)^2}.
\end{equation}
With $x_k=e^{-2k}$ one has $Q_k(1/e)=\big[(k-1)-2(k^2-1)x_k-(k+1)x_k^2\big]/(k^2-1)$ and $S_k(1/e)=\big[(k-1)+(k+1)x_k\big]/(k^2-1)$, which gives the formula of Proposition~C.

\subsection{Second derivation: domain expansion}
We check \eqref{eq:Ak1} by a direct expansion in the physical plane. Let $\widetilde\Omega_\varepsilon=\{r<1+\varepsilon\cos k\theta\}$, whose area is $\pi(1+\varepsilon^2/2)$, and write $\mathrm c=\cos k\theta$. Let
\[
\partial K=\{r=s+\varepsilon b\,\mathrm c+\varepsilon^2g_2(\theta)\},\qquad \int_0^{2\pi}g_2\,d\theta=0,
\]
so that $s$ is the mean radius, and set $\ell=\log s<0$. We look for
\[
u=u_0+\varepsilon u_1+\varepsilon^2u_2,\qquad u_0=\frac{\log r}{\ell},\qquad u_1=(\alpha r^k+\beta r^{-k})\,\mathrm c,
\]
where $u_2$ has a radial part $c_0+d_0\log r$ plus modes that average to zero on circles. Below, $\langle\cdot\rangle$ denotes the mean in $\theta$; note $\langle\mathrm c^2\rangle=\frac12$.

\emph{First order.} The conditions $u=0$ on the outer boundary, $u=1$ on $\partial K$ and $|\nabla u|$ constant on $\partial K$ give
\[
\alpha+\beta=-\frac1\ell,\qquad \alpha s^k+\beta s^{-k}=-\frac{b}{s\ell},\qquad k\,(\alpha s^k-\beta s^{-k})=\frac{b}{s\ell}.
\]
With $B=b/(s\ell)$ and $D(s)=(k-1)s^{-k}+(k+1)s^k$, the solution is
\[
\alpha s^k=-\frac{(k-1)B}{2k},\qquad \beta s^{-k}=-\frac{(k+1)B}{2k},\qquad B=\frac{2k}{\ell\,D(s)},\qquad b=\frac{2ks}{D(s)}.
\]
In particular $\lambda_1=0$.

\emph{Second order, outer boundary.} Expanding $u(1+\varepsilon\mathrm c,\theta)=0$ to order $\varepsilon^2$ and averaging gives
\[
-\frac{1}{4\ell}+\frac{k(\alpha-\beta)}2+c_0=0,\qquad\text{hence}\qquad c_0=\frac1{4\ell}+\frac B4\big[(k-1)s^{-k}-(k+1)s^k\big].
\]

\emph{Second order, inner boundary.} With $\eta=\varepsilon b\mathrm c+\varepsilon^2g_2$ one has $u_0(s+\eta)=1+(\eta/s-\eta^2/(2s^2))/\ell+O(\varepsilon^3)$. The averaged $\varepsilon^2$-part of $u=1$ on $\partial K$ reads
\[
-\frac{b^2}{4s^2\ell}+\frac{bk}{2s}(\alpha s^k-\beta s^{-k})+c_0+d_0\ell=0 .
\]
Since $b^2/(4s^2\ell)=B^2\ell/4$ and $\frac{bk}{2s}(\alpha s^k-\beta s^{-k})=B^2\ell/2$, this gives $d_0=-c_0/\ell-B^2/4$.

\emph{Second order, gradient condition.} On $\partial K$ we use $|\nabla u|^2=(\partial_ru)^2+r^{-2}(\partial_\theta u)^2$. The first-order part of $\partial_ru$ vanishes by the first-order equations, so
\[
\partial_ru=\frac1{s\ell}+\varepsilon^2X+O(\varepsilon^3),\qquad
X=-\frac{g_2}{s^2\ell}+\frac{b^2\mathrm c^2}{s^3\ell}+b\,\mathrm c\,\partial_r^2u_1(s)+\partial_ru_2(s),
\]
and $r^{-2}(\partial_\theta u)^2=\varepsilon^2k^2B^2\sin^2(k\theta)/s^2+O(\varepsilon^3)$, because $\alpha s^k+\beta s^{-k}=-B$. Writing $\lambda=\lambda_0+\varepsilon^2\lambda_2$ with $\lambda_0=-1/(s\ell)>0$ and averaging $\lambda^2=|\nabla u|^2$, we get
\[
\lambda_2=-\langle X\rangle+\frac{k^2B^2}{4s^2\lambda_0}=-\langle X\rangle-\frac{k^2B^2\ell}{4s}.
\]
Using $\partial_r^2u_1(s)=-B(k^2+1)s^{-2}\,\mathrm c$, one finds $\langle X\rangle=-k^2B^2\ell/(2s)+d_0/s$, and therefore
\begin{equation}\label{eq:lambda2s}
\lambda_2(s)=\frac1s\Big[\frac{k^2B^2\ell}{4}+\frac1{4\ell^2}+\frac{B}{4\ell}\big((k-1)s^{-k}-(k+1)s^k\big)+\frac{B^2}{4}\Big].
\end{equation}
At $s=1/e$ ($\ell=-1$, $B=-2k/D$, $D=(k-1)e^k+(k+1)e^{-k}$) this becomes
\[
\lambda_2(1/e)=e\Big[\frac14-\frac k2+\frac{k(k-1)e^k}{D}-\frac{k^2(k^2-1)}{D^2}\Big].
\]
Since $c(t\Omega)=c(\Omega)/t$ and $\Omega_\varepsilon=\widetilde\Omega_\varepsilon/\sqrt{1+\varepsilon^2/2}$, the area normalisation contributes $+e/4$. Hence
\[
A_k=e\Big[\frac{1-k}2+k\Phi_k-k^2\Phi_k(1-\Phi_k)\Big],\qquad \Phi_k=\frac{(k-1)e^{2k}}{(k-1)e^{2k}+k+1}.
\]
Put $x_k=e^{-2k}$ and $D_k'=(k-1)+(k+1)x_k$, so that $\Phi_k=(k-1)/D_k'$. One checks
\[
2D_k'^2\Big[\frac{1-k}2+k\Phi_k-k^2\Phi_k(1-\Phi_k)\Big]=(k-1)^2(k+1)-2(k^2-1)^2x_k-(k-1)(k+1)^2x_k^2 .
\]
The right-hand side is exactly $(k^2-1)\big[(k-1)-2(k^2-1)x_k-(k+1)x_k^2\big]$, so both derivations agree identically.

\subsection{Positivity}
For $k=1$ the factor $k^2-1$ vanishes, so $A_1=0$, as it must be, since translations leave $c$ invariant. For $k\ge2$, $(k+1)x_k=(k+1)e^{-2k}\le3e^{-4}<0.055$, because $k\mapsto(k+1)e^{-2k}$ is decreasing. Hence
\[
(k-1)-2(k^2-1)x_k-(k+1)x_k^2\ \ge\ (k-1)(1-0.11)-0.055\,x_k>0,
\]
and $A_k>0$. As $k\to\infty$, $A_k=e(k+1)/2+O(k^3e^{-2k})$.

By rotational invariance, the second variation is invariant under rotations; on each two-dimensional space spanned by $\cos k\theta$ and $\sin k\theta$ the rotation group acts irreducibly, so the second variation is a multiple of the identity there. Thus it is the diagonal quadratic form
\[
c(\Omega)=e+\varepsilon^2\sum_{k\ge2}A_k\,(a_k^2+b_k^2)+O(\varepsilon^3)\qquad\text{for } \partial\Omega:\ r\propto1+\varepsilon\sum_k(a_k\cos k\theta+b_k\sin k\theta),
\]
which is positive definite modulo translations and, since $A_k\asymp k$, coercive in the $H^{1/2}$ norm of the boundary perturbation.

\begin{remark}\label{rem:rigour}
To turn Proposition~C into a rigorous local minimality statement, one needs: (a) the implicit function theorem for the branch of solutions near the fold, in H\"older spaces, to justify the expansions; (b) a stability argument of Dambrine--Pierre type \cite{DP} to pass from $H^{1/2}$-coercivity of the second variation to strict local minimality in a $C^{2,\alpha}$ neighbourhood (the usual two-norm discrepancy). We do not carry this out here. The numerical values in Table~\ref{tab:second} confirm the formula independently.
\end{remark}

%=====================================================================
\section{A Rellich--Pohozaev identity and the proof of Theorem D}\label{sec:rellich}
%=====================================================================

Throughout this section $\Omega\subset\R^N$, $N\ge2$, is a bounded domain of class $C^{1,\alpha}$, strictly star-shaped with respect to the origin, and $(K,u)$ is a classical solution of \eqref{eq:bernoulli}. Set $D=\Omega\setminus K$. By the Schauder estimates, $u\in C^{1,\alpha}$ up to $\partial\Omega$, and $u\in C^1$ up to $\partial K$ by assumption, so $u\in C^2(D)\cap C^1(\overline D)$. We denote by $\nu$ the outer unit normal to $\Omega$ on $\partial\Omega$ and to $K$ on $\partial K$.

\subsection{Two integral identities}

\begin{lemma}\label{lem:flux}
$\displaystyle\int_{\partial\Omega}|\nabla u|\,d\sigma=\int_{\partial K}|\nabla u|\,d\sigma=\cp(K;\Omega)=\lambda\Per(K)$.
\end{lemma}

\begin{proof}
By the strong maximum principle $0<u<1$ in $D$, so $|\nabla u|=-\partial_\nu u$ on $\partial\Omega$ and $|\nabla u|=\partial_\nu u$ on $\partial K$ (one-sided derivatives from $D$). Since $u$ is harmonic, $\int_{\partial D}\partial_{\nu_D}u\,d\sigma=0$, which gives the first equality. Green's formula $\int_D|\nabla u|^2=\int_{\partial D}u\,\partial_{\nu_D}u\,d\sigma=\int_{\partial K}|\nabla u|\,d\sigma$ and the minimality of $u_K=u$ give the second; the third follows from $|\nabla u|=\lambda$ on $\partial K$.
\end{proof}

\begin{lemma}[Rellich--Pohozaev identity]\label{lem:rellich}
Identity \eqref{eq:rellich} holds.
\end{lemma}

\begin{proof}
Consider the vector field $V=(x\cdot\nabla u)\nabla u-\frac12|\nabla u|^2x$ \cite{Re,Poh}. For harmonic $u$,
\[
\operatorname{div}\big((x\cdot\nabla u)\nabla u\big)=|\nabla u|^2+x\cdot D^2u\,\nabla u,\qquad
\operatorname{div}\big(\tfrac12|\nabla u|^2x\big)=\tfrac N2|\nabla u|^2+x\cdot D^2u\,\nabla u,
\]
so $\operatorname{div}V=\frac{2-N}2|\nabla u|^2$. Since $V\in C(\overline D)\cap C^1(D)$, $\operatorname{div}V\in L^1(D)$ and $\partial D$ is of class $C^1$, the divergence theorem applies (e.g.\ by exhausting $D$ by smooth subdomains).

On each component of $\partial D$, $u$ is constant, so $\nabla u=(\partial_{\nu_D}u)\,\nu_D$ and $V\cdot\nu_D=\frac12(x\cdot\nu_D)|\nabla u|^2$. On $\partial K$, $\nu_D=-\nu$ and $|\nabla u|=\lambda$, while $\int_{\partial K}x\cdot\nu\,d\sigma=N|K|$ by the divergence theorem in $K$. Using Lemma~\ref{lem:flux},
\[
\frac{2-N}{2}\cp(K;\Omega)=\frac12\int_{\partial\Omega}|\nabla u|^2(x\cdot\nu)\,d\sigma-\frac{\lambda^2}{2}N|K|,
\]
which is \eqref{eq:rellich}.
\end{proof}

\subsection{The geometric factor}

\begin{lemma}\label{lem:theta}
$I(\Omega)\ge|\partial\Omega|^2/(N|\Omega|)\ge I(\Omega^*)=|\partial\Omega^*|/R$, where $|\Omega^*|=|B_R|$. Hence $\theta(\Omega)\in(0,1]$, with equality if and only if $\Omega$ is a ball centred at the origin.
\end{lemma}

\begin{proof}
By Cauchy--Schwarz, $|\partial\Omega|^2\le I(\Omega)\int_{\partial\Omega}x\cdot\nu\,d\sigma=I(\Omega)\,N|\Omega|$. The isoperimetric inequality gives $|\partial\Omega|^2/(N|\Omega|)\ge|\partial B_R|^2/(N|B_R|)=|\partial B_R|/R=I(B_R)$. Equality in the isoperimetric inequality forces $\Omega$ to be a ball; equality in Cauchy--Schwarz forces $x\cdot\nu$ to be constant on $\partial\Omega$, so the ball is centred at the origin.
\end{proof}

\subsection{Proof of Theorem D}
\emph{Step 1: a lower bound by the capacity.} By Lemma~\ref{lem:flux} and Cauchy--Schwarz,
\[
\cp(K;\Omega)^2=\Big(\int_{\partial\Omega}|\nabla u|\,d\sigma\Big)^2\le\int_{\partial\Omega}|\nabla u|^2(x\cdot\nu)\,d\sigma\cdot I(\Omega).
\]
Inserting this in \eqref{eq:rellich} gives
\begin{equation}\label{eq:keyineq}
\lambda^2\,N|K|\ \ge\ \Phi_\Omega\big(\cp(K;\Omega)\big),\qquad \Phi_\Omega(t)=\frac{t^2}{I(\Omega)}+(N-2)\,t,
\end{equation}
and $\Phi_\Omega$ is increasing on $[0,\infty)$.

\emph{Step 2: symmetrisation.} Let $R$ and $s<R$ be defined by $|B_R|=|\Omega|$ and $|B_s|=|K|$, and put $C(s)=\cp(\overline{B_s};B_R)$, that is,
\[
C(s)=\frac{2\pi}{\log(R/s)}\ \ (N=2),\qquad C(s)=\frac{(N-2)|\mathbb S^{N-1}|}{s^{2-N}-R^{2-N}}\ \ (N\ge3).
\]
Extend $u$ by $1$ on $K$ and by $0$ outside $\Omega$, and let $u^*$ be its Schwarz symmetrisation. By the P\'olya--Szeg\H{o} inequality \cite{PS,BZ}, $u^*\in W^{1,2}_0(B_R)$ and $\int|\nabla u^*|^2\le\int|\nabla u|^2$. Since $0<u<1$ in $D$, the set $\{u^*=1\}$ is the closed ball of measure $|\{u=1\}|=|K|$, namely $\overline{B_s}$. Thus $u^*$ is admissible for $\cp(\overline{B_s};B_R)$, and
\[
\cp(K;\Omega)\ \ge\ C(s).
\]
By \eqref{eq:keyineq} and the monotonicity of $\Phi_\Omega$,
\begin{equation}\label{eq:afterPS}
\lambda^2\,N|B_s|\ \ge\ \frac{C(s)^2}{I(\Omega)}+(N-2)\,C(s).
\end{equation}

\emph{Step 3: the bound $\sqrt\theta$.} For the concentric radial solution $K=\overline{B_s}$ in $\Omega=B_R$, the gradient $|\nabla u|=C(s)/|\partial B_R|$ is constant on $\partial B_R$ and $x\cdot\nu=R$. Hence both Cauchy--Schwarz steps above are equalities, and \eqref{eq:rellich} reads
\begin{equation}\label{eq:radialid}
\lambda_{\rm rad}(s)^2\,N|B_s|=\frac{C(s)^2}{I(B_R)}+(N-2)\,C(s).
\end{equation}
Since $I(\Omega)=I(B_R)/\theta$ with $\theta=\theta(\Omega)\le1$ (Lemma~\ref{lem:theta}), \eqref{eq:afterPS} gives
\[
\lambda^2N|B_s|\ \ge\ \theta\,\frac{C(s)^2}{I(B_R)}+(N-2)C(s)\ \ge\ \theta\Big[\frac{C(s)^2}{I(B_R)}+(N-2)C(s)\Big]=\theta\,\lambda_{\rm rad}(s)^2N|B_s|.
\]
Hence $\lambda\ge\sqrt\theta\,\lambda_{\rm rad}(s)\ge\sqrt\theta\,c(B_R)$, by \eqref{eq:cball}.

\emph{Step 4: the bound $\kappa_N$ ($N\ge3$).} Dropping the first term in \eqref{eq:afterPS},
\[
\lambda^2\ \ge\ \frac{N-2}{N}\,\frac{C(s)}{|B_s|}=\frac{N-2}{N}\cdot\frac{N(N-2)}{R^2\,t^2\,(1-t^{N-2})},\qquad t=\frac sR .
\]
The function $t\mapsto t^2-t^N$ is maximal on $(0,1)$ at $t^{N-2}=2/N$, where it equals $t^2(N-2)/N$. Hence
\[
\lambda^2\ \ge\ \frac{N-2}{N}\,\Lambda_N(R)^2,\qquad \Lambda_N(R)=\frac{N\,(N/2)^{\frac1{N-2}}}{R}.
\]
Dividing by $c(B_R)=(N-1)^{\frac{N-1}{N-2}}/R$ gives $\lambda\ge\kappa_N\,c(B_R)$ with $\kappa_N$ as in \eqref{eq:kappa}. Since $c(\Omega^*)=c(B_R)$, Theorem~D follows by taking the infimum over all classical solutions. \qed

\begin{remark}\label{rem:DK}
$\Lambda_N(R)$ is the Daners--Kawohl constant of $B_R$ (for $N=2$, $\Lambda_2(R)=2\sqrt e/R$). Comparing $J_\lambda(v)$ with $J_\lambda(0)$ shows that $\Lambda(\Omega)^2=\inf_K\cp(K;\Omega)/|K|$. Dropping only the boundary term in \eqref{eq:rellich}, Step~4 yields more precisely
\[
c(\Omega)\ \ge\ \sqrt{\tfrac{N-2}{N}}\;\Lambda(\Omega)\qquad(N\ge3),
\]
so that, for star-shaped domains, the two constants $c$ and $\Lambda$ are comparable. Combined with $\Lambda(\Omega)\ge\Lambda(\Omega^*)$ \cite{DK}, this gives the bound $\kappa_N$ again.
\end{remark}

\begin{remark}
As $N\to\infty$, $R\,c(B_R)=N-1+\log N+o(1)$ and $R\,\Lambda_N(R)=N+\log(N/2)+o(1)$, so $\kappa_N\to1$. For $\Omega$ close to a ball, $I(\Omega)/I(\Omega^*)-1$ is quadratic in the deformation (it is minimal at the ball), so the bound $\sqrt\theta$ deviates from $1$ only at second order. This is the same order as in Proposition~C, but with the opposite sign: it does not prove local minimality.
\end{remark}

\subsection{Proof of Corollary E}
A bounded convex domain of class $C^{1,\alpha}$ is strictly star-shaped with respect to any interior point $x_0$, since $(x-x_0)\cdot\nu\ge\operatorname{dist}(x_0,\partial\Omega)>0$. By translation invariance of $c$, Theorem~D applies with the origin at $x_0$, and Theorem~A gives the other bound. For the ellipse $x=(a\cos t,b\sin t)$ centred at $x_0=0$, one has $d\sigma=\sqrt{a^2\sin^2t+b^2\cos^2t}\,dt$ and $x\cdot\nu=ab/\sqrt{a^2\sin^2t+b^2\cos^2t}$, so
\[
I=\frac{\pi(a^2+b^2)}{ab},\qquad \theta=\frac{2ab}{a^2+b^2}=\frac{2}{q+q^{-1}},\qquad q=\frac ab .
\]
The inequality $\sqrt\theta>\mu/e$ is equivalent to $q+q^{-1}<2e^2/\mu^2\approx2.39261$, i.e.\ $q<1.8529$. \qed

%=====================================================================
\section{Polarization: an obstruction and what survives}\label{sec:polarization}
%=====================================================================

For a closed half-space $H$ with reflection $\sigma$, the polarization of a set $E$ is
\[
E^H=\big((E\cup\sigma E)\cap H\big)\cup\big((E\cap\sigma E)\cap H^c\big),
\]
and the polarization of a function $v$ is $\max(v,v\circ\sigma)$ on $H$ and $\min(v,v\circ\sigma)$ on $H^c$; see \cite{BrS}. A natural strategy towards \eqref{eq:FR} is to show that $c(\Omega^H)\le c(\Omega)$ and to conclude by approximating $\Omega^*$ by iterated polarizations. If $c$ admits the Beurling-type characterisation
\[
c(\Omega)=\inf_{K\subset\Omega}\ \sup_{\partial K}|\nabla u_K|,
\]
where $u_K$ is the capacitary potential of $K$ in $\Omega$, it would suffice to transport subsolutions, i.e.\ to deduce $\sup_{\partial K^H}|\nabla u_{K^H}|\le\lambda$ in $\Omega^H$ from $\sup_{\partial K}|\nabla u_K|\le\lambda$. For convex domains, such a characterisation is closely related to the subsolution methods of Beurling \cite{Be}, Acker \cite{Ac} and Henrot--Shahgholian \cite{HS2}; our numerical tests (Section~\ref{sec:beurling}) support it. The next lemma shows that the polarized free set itself is never a good competitor. Note that it does not disprove $c(\Omega^H)\le c(\Omega)$: other free sets in $\Omega^H$ may still do better.

\begin{lemma}[corner obstruction]\label{lem:corner}
Assume that $\partial K$ and $\partial(\sigma K)$ are $C^{1,\alpha}$ curves that cross transversally at a point $P\in\operatorname{int}(H^c)\cap\Omega^H$. Then $\sup_{\partial K^H}|\nabla u_{K^H}|=+\infty$.
\end{lemma}

\begin{proof}
Near $P$, $K^H$ coincides with $K\cap\sigma K$, which has at $P$ a convex corner of interior angle $\vartheta\in(0,\pi)$. Hence $\Omega^H\setminus K^H$ has at $P$ an opening angle $\omega=2\pi-\vartheta>\pi$. The function $v=1-u_{K^H}$ is positive and harmonic in $\Omega^H\setminus K^H$ and vanishes on both arcs of $\partial K^H$ through $P$. By the corner asymptotics \cite[Ch.~4--5]{Gr}, after straightening the two arcs by a conformal map, $v=a\,r^{\pi/\omega}\sin(\pi\phi/\omega)+o(r^{\pi/\omega})$, with $a>0$ by the Hopf lemma at corners. Hence $|\nabla v|\sim a(\pi/\omega)r^{\pi/\omega-1}$ along $\partial K^H$, which is unbounded since $\pi/\omega<1$.
\end{proof}

On the side $H$, the union $K\cup\sigma K$ creates reentrant corners of $K^H$, at which the gradient vanishes. The obstruction is therefore entirely located in $H^c$, where the polarized function is a minimum of harmonic functions and is only superharmonic. The critical free sets of non-symmetric convex domains are smooth and convex, and for most $H$ they cross their reflections transversally. What survives is the energy inequality.

\begin{proposition}\label{prop:energy}
For $K\subset\Omega$ set $E_\lambda(K;\Omega)=\cp(K;\Omega)-\lambda^2|K|$. Then, for every $\lambda$,
\[
E_\lambda(K^H;\Omega^H)\ \le\ E_\lambda(K;\Omega).
\]
\end{proposition}

\begin{proof}
Extend $u_K$ by $0$ outside $\Omega$; it is continuous on $\R^N$ and $\{u_K>0\}=\Omega$, $\{u_K\ge1\}=K$. Its polarization $u_K^H$ is continuous, vanishes outside $\{u_K>0\}^H=\Omega^H$, and therefore belongs to $W^{1,2}_0(\Omega^H)$. It satisfies $u_K^H\ge1$ on $\{u_K\ge1\}^H=K^H$, and has the same Dirichlet integral as $u_K$ \cite{BrS}. Hence $\cp(K^H;\Omega^H)\le\cp(K;\Omega)$, while $|K^H|=|K|$.
\end{proof}

By the Hadamard formula $\delta\cp(K;\Omega)=\int_{\partial K}|\nabla u_K|^2V_n\,d\sigma$, solutions of \eqref{eq:bernoulli} are exactly the critical points of $E_\lambda$. For the disc, the elliptic solution is a local minimum and the hyperbolic one a saddle, and they merge into an inflection point at $\lambda=c$. Proposition~\ref{prop:energy} implies that polarization decreases the Daners--Kawohl constant $\Lambda$ (the threshold for $\inf E_\lambda<0$). It also decreases the capacitary profile $C_\Omega(v)=\inf\{\cp(K;\Omega):|K|=v\}$. But $c$ is the threshold of existence of a critical point, not an energy level, and it is not monotone under a decrease of $E_\lambda$.

In the convex class, Cardaliaguet and Tahraoui \cite{CT} studied the profile
\[
F_\Omega(\sigma)=\inf\{\cp(D;\Omega):\ D\Subset\Omega\ \text{convex},\ |D|=\sigma\}.
\] They showed that, at points where $F_\Omega$ is differentiable, convex minimisers are solutions of level $\sqrt{F_\Omega'(\sigma)}$. This suggests the following question.

\begin{question}
Find a characterisation of $c(\Omega)$ in terms of $E_\lambda$ (or of $F_\Omega$) which is monotone under pointwise decrease of the energy. For instance: is
\[
c(\Omega)=\inf\{\lambda:\ E_\lambda(\cdot\,;\Omega)\ \text{has a strict local minimiser}\},
\]
and does a strict local minimiser for $\Omega$ produce one for $\Omega^H$?
\end{question}

%=====================================================================
\section{Numerical experiments}\label{sec:numerics}
%=====================================================================

\subsection{Method}
We solve \eqref{eq:bernoulli} in the plane by the method of fundamental solutions. The potential is represented as $u(x)=a_0+\sum_ja_j\log|x-z_j|$, with source points $z_j$ placed outside $\Omega$ and inside $K$. The Dirichlet conditions $u=0$ on $\partial\Omega$ and $u=1$ on $\partial K$ are imposed in the least-squares sense at $240$--$320$ collocation points. The free boundary is parametrised in polar form, $\rho(\theta)=s+\sum_{k\le M}(a_k\cos k\theta+b_k\sin k\theta)$, $M=12$--$20$. For fixed mean radius $s$, the equation $|\nabla u|=\lambda$ on $\partial K$ is solved for $(a_k,b_k,\lambda)$ by a Levenberg--Marquardt method. The Bernoulli constant is then the minimum of $\lambda(s)$ along the branch. Every domain is normalised to area $\pi$, so that $c(\Omega^*)=e$.

\emph{Validation.} For the unit disc the method returns $c=2.718281828459045$ at $s=0.36787943$, i.e.\ $e$ and $1/e$ to machine precision. Capacities of concentric annuli are reproduced to $10^{-14}$. Residuals of the free boundary condition are below $10^{-5}$ in all reported runs, and usually below $10^{-10}$.

\subsection{Ellipses and Steiner symmetrisation}
Table~\ref{tab:ellipses} shows that $c$ increases with the aspect ratio of ellipses. Values are stable to $10^{-9}$ between $M=12$ and $M=20$. For aspect ratio $3$ the polar parametrisation of $K$ is no longer adequate, and we do not report it. The last column gives the lower bound of Corollary~E.

\begin{table}[ht]
\centering
\begin{tabular}{@{}cccc@{}}
\toprule
$q$ & $c(\Omega)$ & $c(\Omega)/e$ & $\max\{\mu/e,\sqrt\theta\}$\\
\midrule
1 & 2.718282 & 1 & 1\\
1.2 & 2.745282 & 1.009933 & 0.991769\\
1.25 & 2.758675 & 1.014860 & 0.987730\\
1.5 & 2.850462 & 1.048626 & 0.960769\\
2 & 3.095436 & 1.138747 & 0.914279\\
\bottomrule
\end{tabular}
\caption{Ellipses of area $\pi$ and aspect ratio $q$.}
\label{tab:ellipses}
\end{table}

To test the Steiner route we consider $\Omega_f=\{(x,y):|y-f(x)|<0.8\sqrt{1-x^2}\}$, whose Steiner symmetral with respect to $y$ is the ellipse with semi-axes $1$ and $0.8$. Table~\ref{tab:steiner} shows that symmetrisation decreases $c$ in all cases, including a non-convex one.

\begin{table}[ht]
\centering
\begin{tabular}{@{}lcc@{}}
\toprule
Domain & $f(x)$ & $c/e$\\
\midrule
Steiner symmetral (ellipse) & $0$ & 1.014860\\
Mild banana & $0.25x^2$ & 1.024310\\
S-shape & $0.4x^3$ & 1.032890\\
Strong banana (non-convex) & $0.5x^2$ & 1.049881\\
\bottomrule
\end{tabular}
\caption{Steiner symmetrisation; values normalised to area $\pi$.}
\label{tab:steiner}
\end{table}

\subsection{Second variation}
For $\Omega_\varepsilon$ as in Proposition~C, we computed $(c(\Omega_\varepsilon)-e)/\varepsilon^2$ at $\varepsilon=0.03$ and $\varepsilon=0.06$ and Richardson-extrapolated in $\varepsilon^2$, which is legitimate since the expansion is even in $\varepsilon$. The agreement with the closed formula is about $10^{-4}$ (Table~\ref{tab:second}).

\begin{table}[ht]
\centering
\begin{tabular}{@{}ccccc@{}}
\toprule
$k$ & $\varepsilon=0.03$ & $\varepsilon=0.06$ & extrapolated & $A_k$ (Prop.~C)\\
\midrule
2 & 3.25860 & 3.26208 & 3.2574 & 3.257430\\
3 & 5.27815 & 5.28396 & 5.2762 & 5.276246\\
4 & 6.76165 & 6.75119 & 6.7651 & 6.765339\\
5 & 8.13388 & 8.08991 & 8.1485 & 8.149293\\
\bottomrule
\end{tabular}
\caption{Second variation coefficients.}
\label{tab:second}
\end{table}

\subsection{A capacitary Cheeger inequality fails}
Since $\cp(K;\Omega)=\int_{\partial K}|\nabla u_K|\,d\sigma$, one has $\sup_{\partial K}|\nabla u_K|\ge\cp(K;\Omega)/\Per(K)$, with equality for concentric discs. One might hope to prove \eqref{eq:FR} via $\inf_K\cp(K;\Omega)/\Per(K)\ge c(\Omega^*)$. This is false. Oscillating boundaries make $\Per(K)$ arbitrarily large while the capacity stays bounded, so the infimum is $0$ for every $\Omega$. In the unit disc, already $\rho(\theta)=s(1+0.2\cos4\theta)$ gives $\cp/\Per=2.5169<e$, while $\sup|\nabla u_K|=3.9453>e$. The supremum cannot be replaced by a mean.

\subsection{Testing the Beurling characterisation}\label{sec:beurling}
In the unit disc, $129$ random compact sets $K$ (five Fourier modes, mean radius in $[0.2,0.45]$) all satisfy $\sup_{\partial K}|\nabla u_K|\ge3.65>e$. Random perturbations of the critical free set $K_c$, with $40$ samples for each amplitude, always increase the supremum. The smallest increments are $0.12$, $0.64$ and $1.43$ for amplitudes $0.002$, $0.01$ and $0.03$ in the disc, and $0.18$, $0.64$ and $1.67$ for the ellipse $q=2$.

\subsection{Testing Proposition B on exact solutions}
The converse part of the dictionary turns every element of $\mathcal A$ into an exact solution. We sampled $\rho\in[0.15,0.8]$ and $2$ to $5$ modes $c_n$ in \eqref{eq:laurent}, and adjusted $c_1$ to enforce $\Res e^\varphi=0$. We computed $F$ and the areas spectrally on each circle (with $1024$ Fourier modes), and checked univalence (simple, positively oriented, nested boundary curves). Among $313$ univalent random samples, the minimum of $\lambda R/e$ is $1.000032$. We also ran $24$ Nelder--Mead minimisations of $\lambda R/e$ under the univalence constraint, starting from random configurations. All stopped above $1$, the lowest value being $1.0012$ near the disc ($\rho=0.366$). No violation of \eqref{eq:extremal} was found.

\subsection{The Rellich--Pohozaev identity}
For the critical solution in the ellipse $q=1.5$, the two sides of \eqref{eq:rellich} (with $N=2$, i.e.\ $2\lambda^2|K|=\int_{\partial\Omega}|\nabla u|^2x\cdot\nu\,d\sigma$) evaluate to $6.72954776011166$ and $6.72954776011170$. The bound of Theorem~D, $\sqrt\theta=0.960769$, is compared with the computed values of $c/e$ in Table~\ref{tab:ellipses}.

%=====================================================================
\section{Concluding remarks and open problems}\label{sec:open}
%=====================================================================

\begin{enumerate}
\item \emph{Sharp constant in Theorem A.} The constant $0.914$ comes from neglecting $|K|$ and from Jensen's inequality. A quantitative stability estimate, saying that circle means of $|F'|^2$ close to $m^2$ force $F$ to be close to a homothety and hence $|K|\approx\pi m^2\rho^2$, would improve it. Proving \eqref{eq:extremal} would settle \eqref{eq:FR} for convex planar domains.
\item \emph{Higher dimensions.} Improve the constant $\kappa_N$ of Theorem~D, in particular for $N=3$, where $\kappa_3\approx0.650$. The boundary term $\int_{\partial\Omega}|\nabla u|^2x\cdot\nu\,d\sigma$, which is only estimated through $I(\Omega)$, carries the missing information.
\item \emph{Local minimality.} Make Proposition~C rigorous along the lines of Remark~\ref{rem:rigour}. The domain expansion of Section~\ref{sec:second} extends to $N\ge3$ with spherical harmonics, $u_1=(\alpha r^\ell+\beta r^{2-N-\ell})Y_\ell$; it remains to extend it to the $p$-Laplacian.
\item \emph{Polarization.} Answer the Question of Section~\ref{sec:polarization}.
\item \emph{Riemannian setting.} Bernoulli problems on manifolds have been studied by Djit\'e and Seck \cite{DS} (exterior problem, $p$-Laplacian, compact manifolds) and, in the one-phase framework, by Peng, Zhang and Zhu \cite{PZZ}. The argument of \cite{DK} only uses the characterisation of $\Lambda$ as a capacity/volume quotient and the P\'olya--Szeg\H{o} inequality, which hold for geodesic symmetrisation on $\mathbb S^n$ and $\mathbb H^n$. We therefore expect $\Lambda(\Omega)\ge\Lambda(B)$ for geodesic balls $B$ of the same volume in these space forms. The classical constant $c$ is open there as well.
\end{enumerate}

\subsection*{Acknowledgements}
The authors acknowledge the use of Claude (Anthropic), an AI assistant, for symbolic and numerical verifications, for implementing the numerical experiments and for help in drafting and revising the manuscript. All results, proofs and computations were checked by the authors, who take full responsibility for the content of this paper.


\begin{thebibliography}{99}

\bibitem{Ac} A.~Acker, \emph{Uniqueness and monotonicity of solutions for the interior Bernoulli free boundary problem in the convex $n$-dimensional case}, Nonlinear Anal.\ \textbf{13} (1989), 1409--1425.

\bibitem{AC} H.~W.~Alt and L.~A.~Caffarelli, \emph{Existence and regularity for a minimum problem with free boundary}, J.\ Reine Angew.\ Math.\ \textbf{325} (1981), 105--144.

\bibitem{Be} A.~Beurling, \emph{On free-boundary problems for the Laplace equation}, Seminars on Analytic Functions, vol.~1, Institute for Advanced Study, Princeton, 1957, 248--263.

\bibitem{BS} C.~Bianchini and P.~Salani, \emph{Concavity properties for elliptic free boundary problems}, arXiv:0810.4842.

\bibitem{BrS} F.~Brock and A.~Yu.~Solynin, \emph{An approach to symmetrization via polarization}, Trans.\ Amer.\ Math.\ Soc.\ \textbf{352} (2000), 1759--1796.

\bibitem{BZ} J.~E.~Brothers and W.~P.~Ziemer, \emph{Minimal rearrangements of Sobolev functions}, J.\ Reine Angew.\ Math.\ \textbf{384} (1988), 153--179.

\bibitem{CT} P.~Cardaliaguet and R.~Tahraoui, \emph{Some uniqueness results for Bernoulli interior free-boundary problems in convex domains}, Electron.\ J.\ Differential Equations \textbf{2002}, No.~102, 1--16.

\bibitem{DP} M.~Dambrine and M.~Pierre, \emph{About stability of equilibrium shapes}, M2AN Math.\ Model.\ Numer.\ Anal.\ \textbf{34} (2000), 811--834.

\bibitem{DK} D.~Daners and B.~Kawohl, \emph{An isoperimetric inequality related to a Bernoulli problem}, Calc.\ Var.\ Partial Differential Equations \textbf{39} (2010), 547--555.

\bibitem{DS} A.~S.~Djit\'e and D.~Seck, \emph{Shape optimization and Bernoulli free boundary problems for the Riemannian $p$-Laplacian}, arXiv:2609.00571.

\bibitem{FR} M.~Flucher and M.~Rumpf, \emph{Bernoulli's free-boundary problem, qualitative theory and numerical approximation}, J.\ Reine Angew.\ Math.\ \textbf{486} (1997), 165--204.

\bibitem{Gr} P.~Grisvard, \emph{Elliptic Problems in Nonsmooth Domains}, Pitman, Boston, 1985.

\bibitem{HS1} A.~Henrot and H.~Shahgholian, \emph{Existence of classical solutions to a free boundary problem for the $p$-Laplace operator. I. The exterior convex case}, J.\ Reine Angew.\ Math.\ \textbf{521} (2000), 85--97.

\bibitem{HS2} A.~Henrot and H.~Shahgholian, \emph{Existence of classical solutions to a free boundary problem for the $p$-Laplace operator. II. The interior convex case}, Indiana Univ.\ Math.\ J.\ \textbf{49} (2000), 311--323.

\bibitem{LS} I.~Ly and D.~Seck, \emph{Isoperimetric inequality for an interior free boundary problem with $p$-Laplacian operator}, Electron.\ J.\ Differential Equations \textbf{2004}, No.~109, 1--12.

\bibitem{PZZ} Yu.~Peng, H.-C.~Zhang and X.-P.~Zhu, \emph{One-phase free boundary problems on Riemannian complexes}, Tr.\ Mat.\ Inst.\ Steklova \textbf{330} (2025).

\bibitem{Poh} S.~I.~Pohozaev, \emph{On the eigenfunctions of the equation $\Delta u+\lambda f(u)=0$}, Soviet Math.\ Dokl.\ \textbf{6} (1965), 1408--1411.

\bibitem{PS} G.~P\'olya and G.~Szeg\H{o}, \emph{Isoperimetric Inequalities in Mathematical Physics}, Princeton University Press, 1951.

\bibitem{Po} Ch.~Pommerenke, \emph{Boundary Behaviour of Conformal Maps}, Springer, Berlin, 1992.

\bibitem{Re} F.~Rellich, \emph{Darstellung der Eigenwerte von $\Delta u+\lambda u=0$ durch ein Randintegral}, Math.\ Z.\ \textbf{46} (1940), 635--636.

\bibitem{Se16} D.~Seck, \emph{On an isoperimetric inequality and various methods for the Bernoulli's free boundary problems}, S\~ao Paulo J.\ Math.\ Sci.\ \textbf{10} (2016), 36--59.

\end{thebibliography}
\end{document}